\documentclass[nobibstyle,STIX1COL]{WileyNJD-v2}

\usepackage{mathtools}

\articletype{Article Type}%

\received{TODO}
\revised{TODO}
\accepted{TODO}

\DeclareMathOperator{\vecop}{vec}
\DeclarePairedDelimiter{\norm}{\lVert}{\rVert}
\DeclarePairedDelimiter{\abs}{\lvert}{\rvert}
\newcommand{\LHP}{\mathrm{LHP}}
\newcommand{\RHP}{\mathrm{RHP}}

\begin{document}

\title{Iterative and doubling algorithms for Riccati-type matrix equations: a comparative introduction\protect\thanks{
F.~Poloni acknowledges the support of Istituto Nazionale di Alta Matematica (INDAM), and of a PRA (\emph{progetti di ricerca di ateneo}) project of the University of Pisa.
}}

\author[1]{Federico Poloni*}

\authormark{Federico Poloni}

\address[1]{\orgdiv{Dipartimento di Informatica}, \orgname{Universit\`a di Pisa}, \orgaddress{\state{Pisa}, \country{Italy}}}

\corres{*\email{federico.poloni@unipi.it}}

\presentaddress{Dipartimento di Informatica, Largo Pontecorvo, 56127 Pisa, Italy.}

\abstract{
	We review a family of algorithms for Lyapunov- and Riccati-type equations which are all related to each other by the idea of \emph{doubling}: 
	they construct the iterate $Q_k = X_{2^k}$ of another naturally-arising fixed-point iteration $(X_h)$ via a sort of repeated squaring. 
	The equations we consider are Stein equations $X - A^*XA=Q$, Lyapunov equations $A^*X+XA+Q=0$, discrete-time algebraic Riccati equations $X=Q+A^*X(I+GX)^{-1}A$, continuous-time algebraic Riccati equations $Q+A^*X+XA-XGX=0$, palindromic quadratic matrix equations $A+QY+A^*Y^2=0$, and nonlinear matrix equations $X+A^*X^{-1}A=Q$.
	We draw comparisons among these algorithms, highlight the connections between them and to other algorithms such as subspace iteration, and discuss open issues in their theory.
}

\keywords{doubling algorithm, algebraic Riccati equation, control theory, numerical linear algebra}

\jnlcitation{\cname{%
\author{Poloni F.}} (\cyear{TODO}), 
\ctitle{TODO: shouldn't these fields be filled in automatically?}, \cjournal{TODO}, \cvol{TODO}.}

\maketitle

\section{Introduction}

Riccati-type matrix equations are a family of matrix equations that appears very frequently in literature and applications, especially in systems theory. One of the reasons why they are so ubiquitous is that they are equivalent to certain invariant subspace problems; this equivalence connects them to a larger part of numerical linear algebra, and opens up avenues for many solution algorithms. 

Many books (and even more articles) have been written on these equations; among them, we recall the classical monography by Lancaster and Rodman~\cite{LancasterRodman}, a review book edited by Bittanti, Laub and Willems~\cite{BittantiLaubWillems}, various treatises which consider them from different points of view such as~\cite{AbouKandil,bart1,bart2,bimbook,BoyEFB-book,Datta,IonescuOaraWeiss,Meh91-book}, and recently also a book devoted specifically to doubling~\cite{doublingbook}. 

This vast theory can be presented from different angles; in this exposition, we aim to present a selection of topics which differs from that of the other books and treatises. We focus on introducing doubling algorithms with a direct approach, explaining in particular that they arise as `doubling variants' of other more basic iterations, and detailing how they are related to the subspace iteration, to ADI, to cyclic reduction and to Schur complements. We do not treat algorithms and equations with the greatest generality possible, to reduce technicalities; we try to present the proofs only up to a level of detail that makes the results plausible and allows the interested reader to fill the gaps.

The basic idea behind doubling algorithms can be explained through the `model problem' of computing $w_h = M^{2^h}v$ for a certain matrix $M\in\mathbb{C}^{n\times n}$, $v\in\mathbb{C}^{n}$, and $h\in \mathbb{N}$. There are two possible ways to approach this computation:
\begin{enumerate}
	\item[(1)] Compute $v_{k+1} = Mv_k$, for $k=0,1,\dots,2^{h-1}$ starting from $v_0 = v$; then the result is $w_h = v_{2^h}$.
	\item[(2)] Compute $M_{k+1} = (M_k)^2$, for $k=0,1,\dots,h-1$, starting from $M_0 = M$; then the result is $w_h = M_h v$ (repeated squaring).
\end{enumerate}
It is easy to verify that $M_k v = v_{2^k}$ for each $k$. Hence $k$ iterations of (2) correspond to $2^k$ iterations of (1). We say that (2) is a \emph{squaring variant}, or \emph{doubling variant}, of (1). Each of the two versions has its own pros and cons, and in different contexts one or the other may be preferred. If $h$ is moderate and $M$ is large and sparse, one should favor variant (1): sparse matrix-vector products can be computed efficiently, while the matrices $M_k$ would become dense rather quickly, and one would need to compute and store all their $n^2$ entries. On the other hand, if $M$ is a dense matrix of non-trivial size (let us say $n \approx 10^3$ or $10^4$) and $h$ is reasonably large, then variant (2) wins: fewer iterations are needed, and the resulting computations are rich in matrix multiplications and BLAS level-3 operations, hence they can be performed on modern computers even more efficiently than their flop counts suggest. This problem is an oversimplified version, but it captures the spirit of doubling algorithms, and explains perfectly in which cases they work best.

Regarding competing methods: we mention briefly in our exposition Newton-type algorithms, ADI, and Krylov-type algorithms. We do not treat here direct methods, incuding Schur decomposition-based methods~\cite{Lau79-schur,PapLS80,Van81}, methods based on structured QR~\cite{Bye86a,BunM86,Meh88}, on symplectic URV decompositions~\cite{BenMX98,ChuLM07}, and linear matrix inequalities~\cite{BoyEFB-book}. Although these competitors may be among the best methods for dense problems, they do not fit the scope of our exposition and they do not lend themselves to an immediate comparison with the algorithms that we discuss.

The equations that we treat arise mostly from the study of dynamical systems, both in discrete and continuous time. In our exposition, we chose to start from the discrete-time versions: while continuous-time Riccati equations are simpler and more common in literature, it is more natural to start from discrete-time problems in this context. Indeed, when we discuss algorithms for continuous-time problems we shall see that often the first step is a reduction to a discrete-time problem (possibly implicit).

% In more detail, we begin our discussion with linear problems, Stein (section TODO) and Lyapunov equations (section TODO). Then we move on to DAREs and CAREs (chapters~TODO and~TODO respectively), and conclude with a couple of closely related 

In the following, we use the notation $ A \succ B$ (resp.~$A \succeq B$) to mean that $A-B$ is positive definite (resp.~semidefinite) (Loewner order). We use $\rho(M)$ to denote the spectral radius of $M$, the symbol $\LHP = \{z\in \mathbb{C}: \operatorname{Re}(z) < 0\}$ to denote the (open) left half-plane, and $\RHP$ for the (open) right half-plane. We use the notation $\Lambda(M)$ to denote the spectrum of $M$, i.e., the set of its eigenvalues. We use $M^*$ to denote the conjugate transpose, and $M^\top$ to denote the transpose without conjugation, which appears when combining vectorizations and Kronecker products with the identity $\operatorname{vec}(MXN) = (N^\top \otimes M)\operatorname{vec}(X)$~\cite[Sections~1.3.6--1.3.7]{gvl}.

\section{Stein equations} \label{sec:stein}

The simplest matrix equation that we consider is the \emph{Stein equation} (or \emph{discrete-time Lyapunov equation}).
\begin{equation} \label{stein}
	X - A^*XA = Q, \quad Q=Q^* \succeq 0,
\end{equation}
for $A,X,Q\in \mathbb{C}^{n\times n}$. This equation often arises in the study of discrete-time constant-coefficient linear systems
\begin{equation} \label{dlinearsystem}
	x_{k+1} = Ax_k.
\end{equation}
A classical application of Stein equations is the following. If $X$ solves~\eqref{stein}, then by multiplying by $x_k^*$ and $x_k$ on both sides one sees that $V(x) := x^* X x$ is decreasing over the trajectories of~\eqref{dlinearsystem}, i.e., $V(x_{k+1}) \leq V(x_k)$. This fact can be used to prove stability of the dynamical system~\eqref{dlinearsystem}.

\subsection{Solution properties}

The Stein equation~\eqref{stein} is linear, and can be rewritten using Kronecker products as 
\begin{equation} \label{stein-linearized}
	(I_{n^2} - A^\top \otimes A^*) \vecop(X) = \vecop(Q).	
\end{equation}
If $A=UTU^*$ is a Schur factorization of $A$, then we can factor the system matrix as
\begin{align} \label{stein-schur}
	I_{n^2} - M &= I_{n^2} - A^\top \otimes A^* = (\bar{U}  \otimes U) (I_{n^2} - T^\top \otimes T^*) (U^\top \otimes U^*),	& M &= A^\top \otimes A^*,
\end{align}
which is a Schur-like factorization where the middle term is lower triangular. One can tell when $I-M$ is invertible by looking at its diagonal entries: $I-M$ is invertible (and hence~\eqref{stein} is uniquely solvable) if and only if $\lambda_i \overline{\lambda_j} \neq 1$ for each pair of eigenvalues $\lambda_i,\lambda_j$ of $A$. This holds, in particular, when $\rho(A) < 1$. When the latter condition holds, we can apply the Neumann inversion formula 
\begin{equation} \label{neumann}
(I-M)^{-1} = I + M + M^2 + \dots, 
\end{equation}
which gives (after de-vectorization) an expression for the unique solution as an infinite series
\begin{equation} \label{steinsol}
X = \sum_{k=0}^{\infty} (A^*)^k Q A^k.
\end{equation}
It is apparent from~\eqref{steinsol} that $X \succeq 0$. A reverse result holds, but with strict inequalities: if~\eqref{stein} holds with $X\succ 0$ and $Q \succ 0$, then $\rho(A) < 1$ \cite[Exercise~7.10]{Datta}.

\subsection{Algorithms}

As discussed in the introduction, we do not describe here direct algorithms of the Bartels--Stewart family~\cite{BarS72,ChuBS,EptonBS,GardinerBS} (which, essentially, exploit the decomposition~\eqref{stein-schur} to reduce the cost of solving~\eqref{stein-linearized} from $\mathcal{O}(n^6)$ to $\mathcal{O}(n^3)$) even if they are often the best performing ones for dense linear (Stein or Lyapunov) equations. Rather, we present here two iterative algorithms, which we will use to build our way towards algorithms for nonlinear equations.

The Stein equation~\eqref{stein} takes the form of a fixed-point equation; this fact suggests the fixed-point iteration
\begin{align} \label{smith}
	X_0 &= 0, &  X_{k+1} &= Q + A^* X_k A,
\end{align}
known as \emph{Smith method}~\cite{Smi68}. It is easy to see that the $k$th iterate $X_k$ is the partial sum of~\eqref{steinsol} (and~\eqref{neumann}) truncated to $k+1$ terms, thus convergence is monotonic, i.e., $Q = X_0 \preceq X_1 \preceq X_2 \preceq \dots \preceq X$. Moreover, some manipulations give
\begin{align*}
\vecop(X-X_k) &= (I+M+M^2+\dots) \vecop(Q) - (I+M+M^2+\dots+M^k)\vecop(Q) \\&= M^{k+1}(I+M+M^2+\dots) \vecop(Q) = M^{k+1} \vecop(X),
\end{align*}

or, devectorizing,
\begin{equation} \label{stein-error}
	X-X_k = (A^*)^{k+1} X A^{k+1}.	
\end{equation}
This relation~\eqref{stein-error} implies $\norm{X-X_k} = \mathcal{O}(r^{k})$ for each $r > \rho(A)^2$, so convergence is linear when $\rho(A) < 1$, and it typically slows down when $\rho(A) \approx 1$.

A doubling variant comes from splitting the partial sums into two halves. The truncated sums of~\eqref{neumann} to $2^{k+1}$ terms can be computed iteratively using the identity 
\[
I+M+M^2+\dots+M^{2^{k+1}-1} = (I+M+M^2+\dots+M^{2^k-1}) + M^{2^k}(I+M+M^2+\dots+M^{2^k-1}),
\]
 without computing all the intermediate sums. Setting $\vecop Q_k := (I+M+M^2+\dots+M^{2^k-1}) \vecop{Q}$ and $A_k := A^{2^k}$, one gets the iteration
\begin{subequations} \label{squared-smith}
\begin{align}
	A_0 &= A,& A_{k+1} &= A_k^2,\\
	Q_0 &= Q,& Q_{k+1} &= Q_k + A_k^* Q_k A_k.
\end{align}
\end{subequations}
In view of the definitions, we have $Q_k = X_{2^k}$; so this method computes the $2^k$th iterate of the Smith method directly with $\mathcal{O}(k)$ operations, without going through all intermediate ones. Convergence is quadratic: $\norm{X-Q_k} = \mathcal{O}(r^{2^k})$ for each $r > \rho(A)^2$. The method~\eqref{squared-smith} is known as \emph{squared Smith}. It has been used in the context of parallel and high-performance computing~\cite{BenQQ}, and reappeared in recent years, when it has been used for large and sparse equations~\cite{Pen99,Sad12,BenES} in combination with Krylov methods.

% [TODO: section on applications, mentioning discrete-time stability and Gramians?]

\section{Lyapunov equations} \label{sec:lyap}

Lyapunov equations 
\begin{equation} \label{lyap}
A^*X + XA + Q = 0, \quad Q =Q^* \succeq 0	
\end{equation}
are the continuous-time counterpart of Stein equations. They arise from the study of continuous-time constant-coefficient linear systems
\begin{equation} \label{clinearsystem}
	\frac{\mathrm{d}}{\mathrm{d}t} x(t) = Ax(t).
\end{equation}
A classical application is the following. If $X$ solves~\eqref{lyap}, by multiplying on by $x(t)^*$ and $x(t)$ on both sides one sees that $V(x):= x^* X x$ is decreasing over the trajectories of~\eqref{clinearsystem}, i.e., $\frac{\mathrm{d}}{\mathrm{d}t} V(x(t))\leq 0$. This fact can be used to prove stability of the dynamical system~\eqref{clinearsystem}. Today stability is more often proved by computing eigenvalues, but Stein equations~\eqref{stein} and Lyapunov equations~\eqref{lyap} have survived in many other applications in systems and control theory, for instance in model order reduction~\cite{BenBD11,GugA04,Sim-review}, or as the inner step in Newton methods for other equations (see for instance~\eqref{newton-step} in the following).

\subsection{Solution properties}

Using Kronecker products, one can rewrite~\eqref{lyap} as
\begin{equation} \label{lyap-linearization}
	(I_n \otimes A^* + A^\top \otimes I_n)\vecop(X) = -\vecop(Q),
\end{equation}
and a Schur decomposition $A=UTU^*$ produces
\begin{equation} \label{lyap-schur}
	I_n \otimes A^* + A^\top \otimes I_n = (\bar{U}  \otimes U) (I_n \otimes T^* + T^\top \otimes I_n) (U^\top \otimes U^*).
\end{equation}
Again, this is a Schur-like factorization, where the middle term is lower triangular. One can tell when $I_n \otimes A^* + A^\top \otimes I_n$ is invertible by looking at its diagonal entries: that matrix is invertible (and hence~\eqref{lyap} is uniquely solvable) if and only if $\bar{\lambda_i} + \lambda_j \neq 0$ for each pair of eigenvalues $\lambda_i,\lambda_j$ of $A$. This holds, in particular, if the eigenvalues of $A$ all lie in $\LHP = \{z \in \mathbb{C} \colon \operatorname{Re}(z) < 0 \}$. When the latter condition holds, an analogue of~\eqref{steinsol} is
\begin{equation}
	X = \int_0^{\infty} \exp(A^* t) Q \exp(At) \, \mathrm{d}t.	
\end{equation}
Indeed, this integral converges for every choice of $Q$ if and only if the eigenvalues of $A$ all lie in $\LHP$.

Notice the pleasant symmetry with the Stein case: the (discrete) sum turns into a (continuous) integral; the stability condition for discrete-time linear time-invariant dynamical systems $\rho(A) < 1$ turns into the one $\Lambda(A) \subset \LHP$ for continuous-time systems. Perhaps a bit less evident is the equivalence between the condition $\bar{\lambda_i} + \lambda_j \neq 0$ (i.e., no two eigenvalues of $A$ are mapped into each other by reflection with respect to the imaginary axis) and $\lambda_i \overline{\lambda_j} \neq 1$ (i.e., no two eigenvalues of $A$ are mapped into each other by circle inversion with respect to the complex unit circle).

Lyapunov equations can be turned into Stein equations and \emph{vice versa}. Indeed, for a given $\tau \in \mathbb{C}$, \eqref{lyap} is equivalent to
\[
(A^*-\tau I)X(A-\bar{\tau}I) - (A^*+\bar{\tau} I)X(A+\tau I) - 2\operatorname{Re}(\tau)Q = 0,
\]
or, if $A-\bar{\tau} I$ is invertible,
\begin{align} \label{lyap-to-stein}
X - c(A)^* X c(A) &= 2\operatorname{Re}(\tau)(A^*-\tau I)^{-1} Q (A-\bar{\tau}I)^{-1}, & c(A) &= (A+\tau I)(A-\bar{\tau} I)^{-1}=(A-\bar{\tau} I)^{-1}(A+\tau I).
\end{align}
If $\tau \in \RHP$, then the right-hand side is positive semidefinite and~\eqref{lyap-to-stein} is a Stein equation. The stability properties of $c(A)$ can be explicitly related to those of $A$ via the following lemma.

\begin{lemma}[properties of Cayley transforms] \label{lem:cayley}
Let $\tau \in \RHP$. Then,
\begin{enumerate}
	\item[(1)] for $\lambda \in \mathbb{C}$, we have $\abs{c(\lambda)} = \abs*{\frac{\lambda+\tau}{\lambda - \bar{\tau}}}<1$ if and only if $\lambda \in \LHP$;
	\item[(2)] for a matrix $A\in\mathbb{C}^{n\times n}$, we have $\rho(c(A)) < 1$ if and only if $\Lambda(A) \subset \LHP$.
\end{enumerate}
\end{lemma}
A geometric argument to visualize (1) is the following. In the complex plane, $-\tau$ and $\bar{\tau}$ are symmetric with respect to the imaginary axis, with $-\tau$ lying to its left. Thus a point $\lambda \in \mathbb{C}$ is closer to $-\tau$ than to $\bar{\tau}$ if and only if it lies in $\LHP$. Part (2) follows from facts on the behaviour of eigenvalues of a matrix under rational functions~\cite[Proposition~1.7.3]{LancasterRodman}, which we will often use also in the following.

Another important property of the solutions $X$ of Lyapunov and Stein equations is the decay of their singular values in many practical cases. We defer its discussion to the following section, since a proof follows from the properties of certain solution algorithms.

\subsection{Algorithms}

As in the Stein case, one can implement a direct $\mathcal{O}(n^3)$ Bartels-Stewart algorithm~\cite{BarS72} by exploiting the decomposition~\eqref{lyap-schur}: the two outer factors have Kronecker product structure, and the inner factor is lower triangular, allowing for forward substitution. An interesting variant  allows one to compute the Cholesky factor of $X$ directly from the one of $Q$~\cite{Ham82}.

Again, we focus our interest on iterative algorithms. We will assume $\Lambda(A) \subset \LHP$. Then, thanks to Lemma~\ref{lem:cayley}, we have $\rho(c(A)) < 1$, so we can apply the Smith method~\eqref{smith} to~\eqref{lyap-to-stein}. In addition, we can change the value of $\tau$ at each iteration. The resulting algorithm is known as \emph{ADI iteration}~\cite{PeaR55,Wac88}:
\begin{align} \label{adi}
	X_0 &= 0, & X_{k+1} &= Q_k + c_k(A)^* X_k c_k(A),\\
	Q_k &= 2\operatorname{Re}(\tau_k)(A^*-\tau_k I)^{-1} Q (A-\bar{\tau}_k I)^{-1}, &  c_k(A) &= (A+\tau_k I)(A-\bar{\tau}_k I)^{-1}=(A-\bar{\tau}_k I)^{-1}(A+\tau_k I). \nonumber
\end{align}
The sequence of \emph{shifts} $\tau_k \in \RHP$ can be chosen arbitrarily, with the only condition that $\bar{\tau}_k \not \in \Lambda(A)$. By writing a recurrence for the error $E_k = X - X_k$, one sees that
\begin{equation} \label{adi-error}
	E_k = r_{k+1}(A)^* E_0 r_{k+1}(A) = r_{k+1}(A)^* X r_{k+1}(A), \quad r_{k+1}(A) = c_k(A) \dots c_1(A)c_0(A),
\end{equation}
a formula which generalizes~\eqref{stein-error}. When $A$ is normal, the problem of assessing the convergence speed of this iteration can be reduced to a scalar approximation theory problem. Note that
\begin{align*}
\norm{r_k(A)} &= \max_{\lambda \in \Lambda(A)} \abs{r_k(\lambda)},  & \norm{r_k(A)^*} = \norm{r_k(-A^*)^{-1}} = \frac{1}{\min_{\lambda \in \Lambda(A)} \abs{r_k(-\lambda^*)} }.
 \end{align*}
If one knows a region $E \subset \LHP$ that encloses the eigenvalues of $A$, the optimal choice of $r_k$ is the degree-$k$ rational function that minimizes
\begin{equation} \label{zolotarev-objective}
	\frac{\sup_{z\in E} \abs{r_k(z)}}{\inf_{z\in -E^*} \abs{r_k(z)}},	
\end{equation}
i.e., a rational function that is `as large as possible' on $E$ and `as small as possible' on $-E^*$.
Finding this rational function is known as \emph{Zolotarev approximation problem}, and it was solved by its namesake for many choices of $E$, including $E=[a,b] \subseteq \mathbb{R}_+$: this choice of $E$ corresponds to having a symmetric positive definite $A$ for which a lower and upper bound on the spectrum are known. It is known that the optimal ratio~\eqref{zolotarev-objective} decays as $\rho^k$, where $\rho < 1$ is a certain value that depends on $E$, related to its so-called \emph{logarithmic capacity}. See the recent review by Beckermann and Townsend~\cite{BecT19} for more details. Optimal choices for the shifts for a normal $A$ were originally studied by Wachspress~\cite{Wac88,EllW91}. When $A$ is non-normal, a similar bound can be obtained from its eigendecomposition $A = VDV^{-1}$, but it includes its eigenvalue condition number $\kappa(V) = \norm{V} \norm{V}^{-1}$, and thus it is of worse quality.

An important case, both in theory and in practice, is when $Q$ has low rank. One usually writes $Q = C^*C$, where $C\in\mathbb{C}^{p\times n}$ is a short-fat matrix, motivated by a standard notation in control theory. A decomposition $X_k = Z_k Z_k^*$ can be derived from~\eqref{adi}, and reads
\begin{align}
Z_k &= \begin{bmatrix}
    \sqrt{2\operatorname{Re}(\tau_{k-1})} (A^*-\tau_{k-1}I)^{-1}C^*, & c_{k-1}(A)^* Z_{k-1} 
\end{bmatrix} \nonumber\\
& = \left[
    \sqrt{2\operatorname{Re}(\tau_{k-1})} (A^*-\tau_{k-1}I)^{-1}C^*, \sqrt{2\operatorname{Re}(\tau_{k-2})}(A^*-\tau_{k-1}I)^{-1}(A^*+\bar{\tau}_{k-1}I) (A^*-\tau_{k-2}I)^{-1}C^*, \dots, \right. \nonumber \\
& \left. 
\sqrt{2\operatorname{Re}(\tau_{0})}(A^*-\tau_{k-1}I)^{-1}(A^*+\bar{\tau}_{k-1}I) (A^*-\tau_{k-2}I)^{-1}(A^*+\bar{\tau}_{k-2}I) \dotsm (A^*-\tau_0 I)^{-1} C^*
\right]. \label{krylovadi}
\end{align}
Hence $Z_k$ is obtained by concatenating horizontally $k$ terms $V_1,V_2,\dots,V_k$ of size $n\times p$ each. Each of them contains a rational function of $A^*$ of increasing degree multiplied by $C^*$. All the factors in parentheses commute: hence that the factors $V_j$ can be computed with the recurrence
\begin{align}
Z_k &= \begin{bmatrix}
    V_1 & V_2 & \cdots & V_{k}
\end{bmatrix}, & V_1 &=\sqrt{2\operatorname{Re}(\tau_{k-1})} (A^*-\tau_{k-1}I)^{-1}C^*, \nonumber \\&& V_{j+1} &= \frac{\sqrt{2\operatorname{Re}(\tau_{k-j-1})}}{\sqrt{2\operatorname{Re}(\tau_{k-j})}}(A^*-\tau_{k-j-1}I)^{-1}(A^*+\bar{\tau}_{k-j}I)V_j \nonumber \\
&&& =
\frac{\sqrt{2\operatorname{Re}(\tau_{k-j-1})}}{\sqrt{2\operatorname{Re}(\tau_{k-j})}} \left(V_j + (\tau_{k-j-1}+\bar{\tau}_{k-j})(A^*+\tau_{k-j-1}I)^{-1}V_j\right). \label{adi-step}
\end{align}
This version of ADI is known as \emph{low-rank ADI (LR-ADI)} \cite{BenLP08}. After $k$ steps, $X_k = Z_kZ_k^*$, but note that in the intermediate steps $j<k$ the quantity $\begin{bmatrix}
    V_1 & V_2 & \cdots & V_j
\end{bmatrix}\begin{bmatrix}
    V_1 & V_2 & \cdots & V_j
\end{bmatrix}^*$ differs from $X_j$ in~\eqref{adi}. Indeed, in this factorized version the shifts appear in reversed order, starting from $\tau_{k-1}$ and ending with $\tau_0$. Nevertheless, we can use LR-ADI as an iteration in its own right: since we keep adding columns to $Z_k$ at each step, $Z_kZ_k^*$  converges monotonically to $X$. This version is particularly convenient for problems in which $A$ is large and sparse, because in each step we only need to solve $p$ linear systems with a shifted matrix $A^*-\tau I$, and we store in memory only the $n \times kp$ matrix $Z_k$. In contrast, iterations such as~\eqref{squared-smith} are not going to be efficient for problems with a large and sparse $A$, since powers of sparse matrices become dense.

The formula~\eqref{krylovadi} displays the relationship between ADI and certain Krylov methods: since the LR-ADI iterates are constructed by applying rational functions of $A^*$ iteratively to $C^*$, the LR-ADI iterate $Z_k$ lies in the so-called \emph{rational Krylov subspace}~\cite{Ruh84}
\begin{equation} \label{ratksub}
	K_{q,k+1}(A^*, C^*) = \operatorname{span} \{q(A^*)^{-1}p(A^*) C^* : \text{$p$ is a polynomial of degree $ \leq k$}\},	
\end{equation}
constructed with \emph{pole polynomial} $q(z)=(z-\tau_0)(z-\tau_1)\dotsm (z-\tau_{k-1})$. This suggests a different view: what is important is not the form of the ADI iteration, but rather the approximation space $K_{q,k}(A^*, C^*)$ to which its iterates belong. Once one has chosen suitable shifts and computed an orthogonal basis $U_k$ of $K_{q,k+1}(A^*, C^*)$, \eqref{lyap} can be solved via \emph{Galerkin projection}: we seek an iterate $X_k$ of the form $X_k = U_k Y_k U_k^*$, and compute $Y_k$ by solving the projected equation
\begin{align*}
0 = U_k^*(A^*X_k + X_kA + Q)U = (U_k^*A^*U_k) Y_k + Y_k (U_k^* A U_k) + U_k^* Q U_k,
\end{align*}
which is a smaller ($kp\times kp$) Lyapunov equation. 

While the approximation properties of classical Krylov subspaces are related to polynomial approximation, those of rational Krylov subspaces are related to approximation with rational functions, as in the Zolotarev problem mentioned earlier. In many cases, rational approximation has better convergence properties, with an appropriate choice of the shifts. This happens also for Lyapunov equations: algorithms based on rational Krylov subspaces~\eqref{ratksub}~\cite{DruS,DruKS} (including ADI which uses them implicitly) often display better convergence properties than equivalent ones in which $U_k$ is chosen as a basis of a regular Krylov subspace or of an extended Krylov subspace
\begin{equation} \label{extKrylov}
K_{k_1,k_2}(A^*,C^*)=\operatorname{span} \{\ell(A^*)C^* : \text{$\ell$ is a Laurent polynomial of degrees $(k_1,k_2)$}\}.	
\end{equation}
Computing a basis for a rational Krylov subspace~\eqref{ratksub} is more expensive than computing one for an extended Krylov subspace~\eqref{extKrylov}: indeed, the former requires solving linear systems with $A-\tau_k I$ for many values of $k$, while the latter uses multiple linear systems with the same matrix $A$. However, typically, their faster convergence more than compensates for it.  Another remarkable feature is the possibility to use an adaptive procedure based on the residual for shift selection~\cite{DruS}.

See also the analysis in Benner, Li, Truhar~\cite{BenLT}, which shows that Galerkin projection can improve also on the ADI solution.

An important consequence of the convergence of these algorithms is that they can be used to give bounds on the rank of the solution $X$. Since we can find rational functions such that~\eqref{zolotarev-objective} decreases exponentially, the formula~\eqref{adi-error} shows that $X$ can be approximated well with $X_k$, which has rank at most $k \cdot \operatorname{rank}(Q)$ in view of the decomposition~\eqref{krylovadi}. This observation has practical relevance, since in many applications $p$ is very small, and the exponential decay in the singular values of $X$ is very well visible and helps reducing the computational cost.

\subsection{Remarks}

There is vast literature already for linear matrix equations, especially when it comes to large and sparse problems. We refer the reader to the review by Simoncini~\cite{Sim-review} for more details. The literature typically deals with continuous-time Lyapunov equations more often than their discrete-time counterpart; however, Cayley transformations~\eqref{lyap-to-stein} can be used to convert one to the other.

In particular, it follows from our discussion that a step of ADI can be interpreted as transforming the Lyapunov equation~\eqref{lyap} into a Stein equation~\eqref{stein} via a Cayley transform~\eqref{lyap-to-stein} and then applying one step of the Smith iteration~\eqref{smith}. Hence the squared Smith method~\eqref{squared-smith} can be interpreted as a doubling algorithm to construct the ADI iterate $X_{2^k}$ in $k$ iterations only, but with the significant limitation of using \emph{only one shift} $\tau$ in ADI.

It is known that a wise choice of shifts has a major impact on the convergence speed of these algorithms; see e.g. Güttel~\cite{Gut13}. A major challenge for doubling-type algorithms seems incorporating multiple shifts in this framework of repeated squaring. It seems unlikely that one can introduce more than one shift per doubling iteration, but even doing so would be an improvement, allowing one to leverage the theory of rational approximation that underlies ADI and Krylov space methods. 

\section{Discrete-time Riccati equations} \label{sec:dare}

We consider the equation
\begin{align} \label{dare}
X &= Q + A^* X(I+GX)^{-1}A & G=G^*&\succeq 0, & Q=Q^* &\succeq 0, & A,G,Q,X &\in\mathbb{C}^{n\times n},
\end{align}
to be solved for $X = X^* \succeq 0$. This equation is known as \emph{discrete-time algebraic Riccati equation} (DARE), and arises in various problems connected to discrete-time control theory~\cite[Chapter~10]{Datta}. Variants in which $G,Q$ are not necessarily positive semidefinite also exist~\cite{RanT93,Wil71}, but we will not deal with them here to keep our presentation simpler. The non-linear term can appear in various slightly different forms: for instance, if $G = BR^{-1}B^*$ for certain matrices $B\in\mathbb{C}^{n\times m},R\in\mathbb{C}^{m\times m}$, $R=R^* \succ 0$, then one sees with some algebra that
\begin{align} 
	X(I+GX)^{-1} &= (I+XG)^{-1}X =  X - X(I+GX)^{-1}GX \nonumber
	\\&=   X - XBR^{-1/2}(I+R^{-1/2}B^*XBR^{-1/2})^{-1}R^{-1/2}B^*X \nonumber 
	\\&=  X - XB(R+B^*XB)^{-1}B^*X, \label{dare-identities}
\end{align}
and all these forms can be plugged into~\eqref{dare} to obtain a slightly different (but equivalent) equation. In particular, from the versions in the last two rows one sees that $X(I+GX)^{-1}$ is Hermitian, which is not evident at first sight. These identities become clearer if one considers the special case in which $\rho(GX)<1$: in this case, one sees that the expressions in~\eqref{dare-identities} are all different ways to rewrite the sum of the converging series $X-XGX + XGXGX - XGXGXGX + \dots$.

Note that the required inverses exist under our assumptions, because the eigenvalues of $GX$ coincide with those of $G^{1/2}XG^{1/2} \succeq 0$.

\subsection{Solution properties}

For convenience, we assume in the following that $A$ is invertible. The results in this section hold also when it is singular, but to formulate them properly one must deal with matrix pencils, infinite eigenvalues, and generalized invariant subspaces (or \emph{deflating subspaces}), a technical difficulty that we would rather avoid here since it does not add much to our presentation. For a more general pencil-based presentation, see for instance Mehrmann~\cite{meh-cayley}.

For each solution $X$ of the DARE~\eqref{dare}, it holds that
\begin{align} \label{dare-subspace}
\begin{bmatrix}
    A & 0\\
    -Q & I
\end{bmatrix}
\begin{bmatrix}
    I\\X
\end{bmatrix}
 &= 
\begin{bmatrix}
    I & G\\
    0 & A^*
\end{bmatrix}
\begin{bmatrix}
    I\\X
\end{bmatrix}
K, &
K &= (I+GX)^{-1}A.
\end{align}
Equation \eqref{dare-subspace} shows that $\operatorname{Im} \begin{bsmallmatrix}I\\X\end{bsmallmatrix}$ is an \emph{invariant subspace} of
\begin{equation} \label{symplectic}
\mathcal{S} = 
	\begin{bmatrix}
    I & G\\
    0 & A^*
\end{bmatrix}^{-1}
\begin{bmatrix}
    A & 0\\
    -Q & I
\end{bmatrix},
\end{equation}
i.e., $\mathcal{S}$ maps this subspace into itself. In particular, the $n$ eigenvalues (counted with multiplicity) of $K$ are a subset of the $2n$ eigenvalues of $\mathcal{S}$: this can be seen by noticing that the matrix $K$ represents (in a suitable basis) the linear operator $\mathcal{S}$ when restricted to said subspace. Conversely, if one takes a basis matrix $\begin{bsmallmatrix}
    U_1\\U_2
\end{bsmallmatrix}$ for an invariant subspace of $\mathcal{S}$, and if $U_1$ is invertible, then $\begin{bsmallmatrix}
    I\\U_2 U_1^{-1}
\end{bsmallmatrix}$ is another basis matrix, the equality \eqref{dare-subspace} holds, and $X=U_2U_1^{-1}$ is a solution of~\eqref{dare}. Hence,~\eqref{dare} typically has multiple solutions, each associated to a different invariant subspace. However, among them there is a preferred one, which is the one typically sought in applications.

\begin{theorem}\relax\cite[Corollary~13.1.2 and Theorem~13.1.3]{LancasterRodman} \label{thm:daresol}
Assume that $Q \succeq 0$, $G\succeq 0$ and $(A,G)$ is d-stabilizable. Then, \eqref{dare} has a (unique) solution $X_+$ such that
\begin{enumerate}
	\item[(1)] $X_+ = X_+^* \succeq 0$;
	\item[(2)] $X_+ \succeq X$ for any other Hermitian solution $X$;
	\item[(3)] $\rho\left((I+GX_+)^{-1}A\right) \leq 1$.
\end{enumerate}
If, in addition, $(Q,A)$ is d-detectable, then $\rho\left((I+GX_+)^{-1}A\right) < 1$.
\end{theorem}
The hypotheses involve two classical definitions from control theory~\cite{Datta}: \emph{d-stabilizable} (resp. \emph{d-detectable}) means that all Jordan chains of $A$ (resp. $A^*$) that are associated to eigenvalues \emph{outside} the set $\{\abs{\lambda}<1\}$ are contained in the maximal (block) Krylov subspace $\operatorname{span}(B, AB, A^2B, \dots)$ (resp. $\operatorname{span}(C^*, A^*C^*, (A^*)^2C^*, \dots)$). We do not discuss further these hypotheses nor the theorem, which is not obvious to prove; we refer the reader to Lancaster and Rodman~\cite{LancasterRodman} for details, and we just mention that these hypotheses are typically satisfied in control theory applications. This solution $X_+$ is often called \emph{stabilizing} (because of property 3) or \emph{maximal} (because of property 2).

Various properties of the matrix $\mathcal{S}$ in~\eqref{symplectic} follow from the fact that it belongs to a certain class of structured matrices. Let $J = \begin{bsmallmatrix}
    0 & I_n\\
    -I_n & 0
\end{bsmallmatrix}\in\mathbb{C}^{2n\times 2n}$. A matrix $M\in\mathbb{C}^{2n\times 2n}$ is called \emph{symplectic} if $M^*JM = J$, i.e., if it is unitary for the non-standard scalar product associated to $J$. The following properties hold.

\begin{lemma} \label{symplemma}
\begin{enumerate}
	\item[(1)] A matrix in the form~\eqref{symplectic} is symplectic if and only if $G=G^*, Q=Q^*$, and the two blocks called $A,A^*$ in~\eqref{symplectic} are one the conjugate transpose of the other.
	\item[(2)] If $\lambda$ is an eigenvalue of a symplectic matrix with right eigenvector $v$, then $\overline{\lambda}^{-1}$ is an eigenvalue of the same matrix with left eigenvector $v^*J$.
	\item[(3)] Under the hypotheses of Theorem~\ref{thm:daresol} (including the d-detectability one in the end), then the $2n$ eigenvalues of $\mathcal{S}$ are (counting multiplicities) the $n$ eigenvalues $\lambda_1,\lambda_2,\dots,\lambda_n$ of $(I+GX_+)^{-1}A$ inside the unit circle, and the $n$ eigenvalues $\overline{\lambda_i}^{-1}$, $i=1,2,\dots,n$ outside the unit circle. In particular, $\begin{bsmallmatrix}
	    I\\X_+
	\end{bsmallmatrix}$ spans the unique invariant subspace of~$\mathcal{S}$ of dimension $n$ all of whose associated eigenvalues lie in the unit circle. 
\end{enumerate}
\end{lemma}
Parts 1 and 2 are easy to verify from the form~\eqref{symplectic} and the definition of symplectic matrix, respectively. To prove Part 3, plug $X_+$ into~\eqref{dare-subspace} and notice that $K$ has $n$ eigenvalues $\lambda_1,\lambda_2,\dots,\lambda_n$ inside the unit circle; these are also eigenvalues of $\mathcal{S}$. By Part~2, all other eigenvalues lie outside the unit circle.

\subsection{Algorithms}

The shape of~\eqref{dare} suggests the iteration
\begin{align} \label{dare-iterative}
	X_{k+1} &= Q + A^* X_k (I+GX_k)^{-1}A, & X_0 &= 0.
\end{align}
This iteration can be rewritten in a form analogous to~\eqref{dare-subspace}:
\begin{align} \label{dare-subspace-iteration}
\begin{bmatrix}
    A & 0\\
    -Q & I
\end{bmatrix}
\begin{bmatrix}
    I\\X_{k+1}
\end{bmatrix}
 &= 
\begin{bmatrix}
    I & G\\
    0 & A^*
\end{bmatrix}
\begin{bmatrix}
    I\\X_k
\end{bmatrix}
K_k, &
K_k &= (I+GX_k)^{-1}A.
\end{align}
Equivalently, one can write it as
\begin{align} \label{dare-subspace-iteration-form2}
\begin{bmatrix}
    U_{1k}\\
    U_{2k}
\end{bmatrix} &= \mathcal{S}^{-1}\begin{bmatrix}
    I\\X_k
\end{bmatrix}, & \begin{bmatrix}
    I\\ X_{k+1}
\end{bmatrix} = \begin{bmatrix}
    U_{1k}\\
    U_{2k}
\end{bmatrix}(U_{1k})^{-1}.
\end{align}
This form highlights a connection with (inverse) subspace iteration (or orthogonal iteration), a classical generalization of the (inverse) power method to find multiple eigenvalues~\cite{Wat-book}. Indeed, we start from the $2n\times n$ matrix $\begin{bsmallmatrix}
    I\\X_0
\end{bsmallmatrix} = \begin{bsmallmatrix}
    I\\0
\end{bsmallmatrix}$, and at each step we first multiply it by $\mathcal{S}^{-1}$, and then we normalize the result by imposing that the first block is $I$. In inverse subspace iteration, we would make the same multiplication, but then we would normalize the result by taking the $Q$ factor of its QR factorization, instead.

It follows from classical convergence results for the subspace iteration (see e.g. Watkins~\cite[Section~5.1]{Wat-book}) that \eqref{dare-subspace-iteration-form2}~converges to the invariant subspace associated to the $n$ largest eigenvalues (in modulus) of $\mathcal{S}^{-1}$, i.e., the $n$ smallest eigenvalues of $\mathcal{S}$. 
In view of Part~3 of Lemma~\ref{symplemma}, this subspace is precisely $ \operatorname{Im}\begin{bsmallmatrix}
    I\\X_+
\end{bsmallmatrix}$. Note that this unusual normalization is not problematic, since at each step of the iteration (and in the limit) the subspace does admit a basis in which the first $n$ rows form an identity matrix. This argument shows the convergence of~\eqref{dare-iterative} to the maximal solution, under the d-detectability condition mentioned in Theorem~\ref{thm:daresol}, which ensures that there are no eigenvalues on the unit circle.

How would one construct a `squaring' variant of this method? Note that that $\begin{bsmallmatrix}
    U_{1k}\\
    U_{2k}
\end{bsmallmatrix} = \mathcal{S}^{-k} \begin{bsmallmatrix}
    I\\0
\end{bsmallmatrix}$; hence one can think of computing $\mathcal{S}^{-2^k}$ by iterated squaring to obtain $X_{2^k}$ in $k$ steps. However, this idea would be problematic numerically, because it amounts to delaying the normalization in subspace iteration until the very last step. The key to solve this issue is using the LU-like decomposition obtained from~\eqref{symplectic}
\[
\mathcal{S}^{-1} =  \begin{bmatrix}
    A & 0\\
    -Q & I
\end{bmatrix}^{-1}
\begin{bmatrix}
    I & G\\
    0 & A^*
\end{bmatrix}.
\]
We seek an analogous decomposition for the powers of $\mathcal{S}^{-1}$, i.e.,
\begin{equation} \label{ssffact}
	\mathcal{S}^{-2^k} =  \begin{bmatrix}
    A_k & 0\\
    -Q_k & I
\end{bmatrix}^{-1}
\begin{bmatrix}
    I & G_k\\
    0 & A_k^*
\end{bmatrix}.
\end{equation}
The following result shows how to compute this factorization with just one matrix inversion.
\begin{lemma}~\cite{PolR} \label{lem:bmf}
Let $M_1,M_2,N_1,N_2\in\mathbb{C}^{2n\times n}$. The factorization
\begin{align} \label{bmf-factorization}
\begin{bmatrix}
    M_1 & M_2
\end{bmatrix}^{-1} \begin{bmatrix}
    N_1 & N_2
\end{bmatrix} &= \begin{bmatrix}
    A_{11} & 0\\
    A_{21} & I_n
\end{bmatrix}^{-1}
\begin{bmatrix}
    I_n & A_{21}\\
    0 & A_{22}
\end{bmatrix},
&
A_{11},A_{12},A_{21},A_{22} &\in \mathbb{C}^{n\times n}
\end{align}
exists if and only if $\begin{bmatrix}
    N_1 & M_2
\end{bmatrix}$ is invertible, and in that case its blocks $A_{ij}$ are given by
\[
\begin{bmatrix}
    A_{11} & A_{12}\\
    A_{21} & A_{22}
\end{bmatrix}
= 
\begin{bmatrix}
    N_1 & M_2
\end{bmatrix}^{-1}
\begin{bmatrix}
    M_1 & N_2
\end{bmatrix}.
\]
\end{lemma}
A proof follows from noticing that the factorization~\eqref{bmf-factorization} is equivalent to the existence of a matrix $K\in\mathbb{C}^{2n\times 2n}$ such that
\[
K
\begin{bmatrix}
    M_1 & M_2 & N_1 & N_2
\end{bmatrix} = 
\begin{bmatrix}
    A_{11} & 0 & I_n & A_{12}\\
    A_{21} & I_n & 0 & A_{22}
\end{bmatrix},
\]
and rearranging block columns in this expression.

One can apply Lemma~\ref{lem:bmf} (with $[M_1\, M_2] = I$ and $[N_1\, N_2] = \begin{bsmallmatrix}
    I & G_k\\
    0 & A_k^*
\end{bsmallmatrix}
\begin{bsmallmatrix}
    A_k & 0\\
    -Q_k & I
\end{bsmallmatrix}^{-1}$ ) to find a factorization of the term in parentheses in 
\begin{equation} \label{sdaderiv}
\mathcal{S}^{-2^{k+1}} = \mathcal{S}^{-2^k}\mathcal{S}^{-2^k} = \begin{bmatrix}
    A_k & 0\\
    -Q_k & I
\end{bmatrix}^{-1}
\left(
\begin{bmatrix}
    I & G_k\\
    0 & A_k^*
\end{bmatrix}
\begin{bmatrix}
    A_k & 0\\
    -Q_k & I
\end{bmatrix}^{-1}
\right)
\begin{bmatrix}
    I & G_k\\
    0 & A_k^*
\end{bmatrix},
\end{equation}
and use it to construct a decomposition~\eqref{ssffact} of $\mathcal{S}^{-2^{k+1}}$ starting from that of $\mathcal{S}^{-2^k}$. The fact that the involved matrices are symplectic can be used to prove that the relations $A_{11}=A_{22}^*$, $A_{21}=A_{21}^*$, $A_{12}=A_{12}^*$ will hold for the computed coefficients. We omit the details of this computation; what matters are the resulting formulas
\begin{subequations} \label{sda-formulas}
\begin{align}
A_{k+1} &= A_k(I+G_k Q_k)^{-1}A_k,\\
G_{k+1} &= G_k + A_k G_k(I+Q_kG_k)^{-1}A_k^*,\\
Q_{k+1} &= Q_k + A_k^*(I+Q_kG_k)^{-1}Q_kA_k, \label{sda-formulasQ}
\end{align}
with $A_0 = A, Q_0 = Q, G_0 = G$.
\end{subequations}
These formulas are all we need to formulate a `squaring' version of~\eqref{dare-iterative}: for each $k$ it holds that
\[
\mathcal{S}^{-2^k} \begin{bmatrix}
    I_n\\0
\end{bmatrix} = \begin{bmatrix}
    I\\
    Q_k
\end{bmatrix}A_k^{-1},
\]
hence $Q_k = X_{2^k}$, the $2^k$th iterate of~\eqref{dare-iterative}. It is not difficult to show by induction that $0 \preceq Q_0 \preceq Q_1 \preceq \dots \leq Q_k \preceq \dots$, and we have already argued above that $Q_k = X_{2^k} \to X_+$. In view of the interpretation as subspace iteration, the convergence speed of~\eqref{dare-iterative} is linear and proportional to the ratio between the absolute values of the $(n+1)$st and $n$th eigenvalue of $\mathcal{S}$, i.e., between $\sigma := \rho((I+GX_+)A) < 1$ and its inverse $\sigma^{-1}$. The convergence speed of its doubling variant~\eqref{sda-formulas} is then quadratic with the same ratio~\cite{doublingbook}.

The iteration~\eqref{sda-formulas}, which goes under the name of \emph{structure-preserving doubling algorithm}, has been used to solve DAREs and related equations by various authors, starting from Chu, Fan, Lin and Wang~\cite{chu-dare}, but it also appears much earlier: for instance, Anderson~\cite{And78} gave it an explicit system-theoretical meaning as constructing an equivalent system with the same DARE solution. The reader may find in the literature slightly different versions of~\eqref{sda-formulas}, which are equivalent to them thanks to the identities~\eqref{dare-identities}.

More general versions of the factorization~\eqref{ssffact} and of the iteration~\eqref{sda-formulas}, which guarantee existence and boundedness of the iterates under much weaker conditions, have been explored by Mehrmann and~Poloni~\cite{MehP12}. Kuo, Lin and Shieh~\cite{KuoLS} studied the theoretical properties of the factorization~\eqref{ssffact} for general powers $\mathcal{S}^t$, $t\in\mathbb{R}$, drawing a parallel with the so-called \emph{Toda flow} for the QR algorithm.

The limit of the monotonic sequence $0 \preceq G_0 \preceq G_1 \preceq G_2 \preceq \dots$ also has a meaning: it is the maximal solution $Y_+$ of the so-called \emph{dual equation}
\begin{equation} \label{dare-dual}
	Y = G + AY(I+QY)^{-1}A^*,
\end{equation}
which is obtained swapping $Q$ with $G$ and $A$ with $A^*$ in~\eqref{dare}. Indeed, SDA for the DARE~\eqref{dare-dual} is obtained by
swapping $Q$ with $G$ and $A$ with $A^*$ in~\eqref{sda-formulas}, but this transformation leaves the formulas unchanged. The dual equation~\eqref{dare-dual} appears sometimes in applications together with~\eqref{dare}. From the point of view of linear algebra, the most interesting feature of its solution $Y_+$ is that $\begin{bsmallmatrix}
    -Y_+\\
    I
\end{bsmallmatrix}$ is a basis matrix for the invariant subspace associated to the other eigenvalues of $\mathcal{S}$, those outside the unit circle. Indeed, \eqref{ssffact}~gives
\[
\mathcal{S}^{2^k} \begin{bmatrix}
    0\\I
\end{bmatrix} = \begin{bmatrix}
    -G_k\\I_n
\end{bmatrix}A_k^{-*},
\]
so $\begin{bsmallmatrix}
    -Y_+\\
    I
\end{bsmallmatrix}$ is the limit of subspace iteration applied to $\mathcal{S}$ instead of $\mathcal{S}^{-1}$, with initial value $\begin{bsmallmatrix}
    0\\I
\end{bsmallmatrix}$.
In particular, putting all pieces together, the following \emph{Wiener-Hopf factorization} holds
\begin{equation} \label{WHfact}
\mathcal{S} = \begin{bmatrix}
    -Y_+ & I\\
    I & X_+
\end{bmatrix}
\begin{bmatrix}
    \left((I+QY_+)^{-1}A^*\right)^{-1}& 0\\
    0& (I+GX_+)^{-1}A
\end{bmatrix}
\begin{bmatrix}
    -Y_+ & I\\
    I & X_+
\end{bmatrix}^{-1}.	
\end{equation}
This factorization relates explicitly the solutions $X_+, Y_+$ to a block diagonalization of $\mathcal{S}$.

An interesting limit case is the one when only the first part of Theorem~\ref{thm:daresol} holds, $(Q,A)$ is not d-detectable, and the solution $X_+$ exists but $\rho((I+GX_+)A) = 1$. In this case, $\mathcal{S}$ has eigenvalues on the unit circle, and it can be proved that all its Jordan blocks relative to these eigenvalues have even size: one can use a result in Lancaster and Rodman\cite[Theorem~12.2.3]{LancasterRodman}, after taking a factorization $G=BR^{-1}B^*$ with $R\succ 0$ and using another result in the same book~\cite[Theorem~12.2.1]{LancasterRodman} to show that the hypothesis $\Psi(\eta)\succ 0$ holds. 

It turns out that in this case the two iterations still converge, although~\eqref{dare-iterative} becomes sublinear and~\eqref{sda-formulas} becomes linear with rate $1/2$. This is shown by Chiang, Chu, Guo, Huang, Lin and Xu~\cite{Chi09-critical-doubling}; the reader can recognize that the key step there is the study of the subspace iteration in presence of Jordan blocks of even multiplicity.

Note that the case in which the assumptions $Q\succeq 0, G \succeq 0$ do not hold is trickier, because there are examples where~\eqref{dare} does not have a stabilizing solution and $\mathcal{S}$ has Jordan blocks of odd size with eigenvalues on the unit circle: an explicit example is
\begin{align} \label{example-dare-eigenvalues}
A &= \begin{bmatrix}
    1 & 3\\ 0 & 1
\end{bmatrix}, & G &= \begin{bmatrix}
    1 & 1\\
    1 & 1
\end{bmatrix}, &
Q &= \begin{bmatrix}
    1 & 0\\
    0 & -10
\end{bmatrix},
\end{align}
which produces a matrix $\mathcal{S}$ with two simple eigenvalues (Jordan blocks of size $1$) $\lambda_{\pm} \approx 0.598 \pm 0.801i$ with $\abs{\lambda}=1$. Surprisingly, eigenvalues on the unit circle are a generic phenomenon for symplectic matrices, which is preserved under perturbations: a small perturbation of the matrices in~\eqref{example-dare-eigenvalues} will produce a perturbed $\tilde{\mathcal{S}}$ with two simple eigenvalues $\tilde{\lambda}_{\pm}$ that satisfy exactly $\abs{\lambda}=1$, because otherwise Part~2 of Lemma~\ref{symplemma} would be violated.

\section{Continuous-time Riccati equations} \label{sec:care}

We consider the equation 
\begin{align} \label{care}
Q + A^*X + XA - XGX &= 0, & G &= G^* \succeq 0, & Q &= Q^* \succeq 0, & A,G,Q,X &\in \mathbb{C}^{n\times n},
\end{align}
to be solved for $X = X^* \succeq 0$. This equation is known as \emph{continuous-time algebraic Riccati equation} (CARE), and arises in various problems connected to continuous-time control theory~\cite[Chapter~10]{Datta}. Despite the very different form, this equation is a natural analogue of the DARE~\eqref{dare}, exactly like Stein and Lyapunov equations are related to each other.

\subsection{Solution properties}

For each solution $X$ of the CARE, it holds
\begin{align} \label{care-subspace}
\begin{bmatrix}
    A & -G\\
    -Q & -A^*
\end{bmatrix}
\begin{bmatrix}
    I\\X
\end{bmatrix}
&=
\begin{bmatrix}
    I\\X
\end{bmatrix}
M,
&
M &= A-GX.
\end{align}
Hence, $\begin{bsmallmatrix}
    I\\X
\end{bsmallmatrix}$ is an invariant subspace of
\begin{equation} \label{hamiltonian}
\mathcal{H} = \begin{bmatrix}
    A & -G\\
    -Q & -A^*
\end{bmatrix}.
\end{equation}

Like in the discrete-time case, this relation implies that the $n$ eigenvalues of $M$ are a subset of those of $\mathcal{H}$; moreover, we can construct a solution $X= U_2 U_1^{-1}$ to~\eqref{care} from an invariant subspace $\operatorname{Im}\begin{bmatrix}
    U_1\\ U_2
\end{bmatrix}$, whenever $U_1$ is invertible. Among all solutions, there is a preferred one.

\begin{theorem}\relax\cite[Theorems~7.9.1, 9.1.2 and~9.1.5]{LancasterRodman} \label{thm:care-maxsol}
Assume that $Q \succeq 0$, $G\succeq 0$, and $(A,G)$ is c-stabilizable. Then, \eqref{care}~has a (unique) solution $X_+$ such that
\begin{enumerate}
	\item[(1)] $X_+ = X_+^* \succeq 0$;
	\item[(2)] $X_+ \succeq X$ for any other Hermitian solution $X$;
	\item[(3)] $\Lambda(A-GX_+) \subset \overline{\LHP}$.
\end{enumerate}
If, in addition, $(Q,A)$ is c-detectable, then $\Lambda(A-GX_+) \subset \LHP$.
\end{theorem}
$\emph{C-stabilizable}$ and $\emph{c-detectable}$ are defined analogously to their discrete-time counterparts, with the only difference that the domain $\{\abs{\lambda}<1\}$ is replaced by the left half-plane $\LHP$. Again, we do not comment on this theorem, whose proof is not obvious, and refer the reader to Lancaster and Rodman~\cite{LancasterRodman}.

Exactly as in the discrete-time case, various interesting properties of the matrix $\mathcal{H}$ in~\eqref{hamiltonian} follow from the fact that it belongs to a certain class of structured matrices. A matrix $M \in \mathbb{C}^{2n\times 2n}$ is called \emph{Hamiltonian} if $-M^*J = JM$, i.e., if it is skew-self-adjoint with respect to the non-standard scalar product induced by $J$. The following result holds.
\begin{lemma}
\begin{enumerate}
	\item[(1)] A matrix in the form~\eqref{hamiltonian} is Hamiltonian if and only if $G=G^*,Q=Q^*$, and the two matrices called $A,A^*$ in~\eqref{hamiltonian} are one the conjugate transpose of the other.
	\item[(2)] If $\lambda$ is an eigenvalue of a Hamiltonian matrix with right eigenvector $v$, then $-\overline{\lambda}$ is an eigenvalue of the same matrix with left eigenvector $v^*J$.
	\item[(3)] If the hypotheses of Theorem~\ref{thm:care-maxsol} hold (including the c-detectability one), then the $2n$ eigenvalues of $\mathcal{H}$ are (counting multiplicities) the $n$ eigenvalues $\lambda_1,\dots,\lambda_n$ of $A-GX_+$ in the left half-plane, and the $n$ eigenvalues $-\overline{\lambda_i}$, $i=1,\dots,n$ in the right half-plane. In particular, $\begin{bsmallmatrix}
	    I\\X_+
	\end{bsmallmatrix}$ spans the unique invariant subspace of~$\mathcal{H}$ of dimension $n$ all of whose associated eigenvalues lie in the left half-plane. 
\end{enumerate}
\end{lemma}
Parts~1 and~2 are easy to verify from the block decomposition~\eqref{hamiltonian} and the definition of Hamiltonian matrix. To prove Part 3, plug $X_+$ into~\eqref{dare-subspace} and notice that $M$ has $n$ eigenvalues $\lambda_1,\lambda_2,\dots,\lambda_n$ in the left half-plane; these are also eigenvalues of $\mathcal{S}$. By Part~2, all other eigenvalues lie in the right half-plane.

The similarities between~\eqref{care-subspace} and~\eqref{dare-subspace} suggest that CAREs can be turned into DAREs (and \emph{vice versa}) by converting the two associated invariant subspace problems; the ingredient to turn one into the other is the Cayley transform.
\begin{lemma} \label{lem:caretodare}
Let $A,G=G^*, Q=Q^*$ be given, and take $\tau > 0$. Set 
\begin{equation} \label{caretodare}
	\begin{bmatrix}
    A_d & G_d\\
    -Q_d & A_d^*
\end{bmatrix}
=
\begin{bmatrix}
    A-\tau I & -G\\
    Q & A^* -\tau I
\end{bmatrix}^{-1}
\begin{bmatrix}
    A+\tau I & -G\\
    Q & A^* + \tau I
\end{bmatrix}
= I + 2\tau \begin{bmatrix}
    A - \tau I & -G\\
    Q & A^* - \tau I
\end{bmatrix}^{-1}.
\end{equation}
Assume that the inverse exists, and that $A_d$ is invertible. Then, the DARE with coefficients $A_d,G_d,Q_d$ has the same solutions as the CARE with coefficients $A,G,Q$ (and, in particular, the same maximal / stabilizing solution).
\end{lemma}
These formulas~\eqref{caretodare} follow from constructing $\mathcal{S}: = c(\mathcal{H}) = (\mathcal{H}-\tau I)^{-1}(\mathcal{H}+\tau I)$, and then applying Lemma~\ref{lem:bmf} to construct a factorization
\[
\mathcal{S} = \begin{bmatrix}
    I & G_d\\
    0 & A_d^*
\end{bmatrix}^{-1}
\begin{bmatrix}
    A_d & 0\\
    -Q_d & I
\end{bmatrix}.
\]
The matrix $\mathcal{S}$ that we have constructed has the same invariant subspaces as $\mathcal{H}$ because $c(\cdot)$ is an invertible rational function: indeed, from \eqref{care-subspace}, it follows that
\begin{align*}
\mathcal{S}\begin{bmatrix}
    I\\X
\end{bmatrix} = c(\mathcal{H})
\begin{bmatrix}
    I\\X
\end{bmatrix}
&=
\begin{bmatrix}
    I\\X
\end{bmatrix}
c(M),
&
M &= A-GX.
\end{align*}
This relation coincides with~\eqref{dare-subspace}, and shows that a solution $X$ of the CARE is also a solution of the DARE constructed with~\eqref{caretodare}.
Thanks to Lemma~\eqref{lem:cayley}, $M$ has all its eigenvalues in~$\LHP$ if and only if $c(M)$ has all its eigenvalues inside the unit circle, so the stabilizing property of the solution is preserved.

Methods to transform DAREs into CAREs and vice versa based on the Cayley transform appear frequently in the literature starting from the 1960s; see for instance Mehrmann~\cite{meh-cayley}, a paper which explores these transformations and mentions the presence of many ``folklore results'' based on the Cayley transforms, relating the properties of the two associated equations.

Even if we restrict ourselves to the assumption that $A_d$ is invertible when treating the DARE, it is important to remark that Lemma~\ref{lem:caretodare} does not generalize completely to the case when $A_d$ is singular~\cite[Section~6]{meh-cayley}. By considering the poles of $c(\mathcal{H})$ as a function of $\tau$, one sees that $A_d$ is singular if and only if $\tau \in \Lambda(\mathcal{H})$. When this happens, even if $\mathcal{S}$ `exists' in a suitable sense as an equivalent matrix pencil, an invariant subspace of $\mathcal{H}$ for which $\tau \in \Lambda(M)$ cannot be converted to the form~\eqref{dare-subspace}, but only to the subtly weaker form
\begin{align} 
\begin{bmatrix}
    A_d & 0\\
    -Q_d & I
\end{bmatrix}
\begin{bmatrix}
    I\\X
\end{bmatrix}(M-\tau I)
 &= 
\begin{bmatrix}
    I & G_d\\
    0 & A_d^*
\end{bmatrix}
\begin{bmatrix}
    I\\X
\end{bmatrix}
(M+\tau I), &
M &= A-GX.
\end{align}
with an additional singular matrix $M-\tau I$ in the left-hand side. Thus we cannot write the equality~\eqref{dare-subspace}, which identifies $X$ as a solution of the DARE: hence the DARE has fewer solutions than the CARE. The stabilizing solution is always preserved by this transformation, though, because $\Lambda(M) \subset \LHP$ cannot contain $\tau > 0$.

\subsection{Algorithms}

In view of the relation between DAREs and CAREs that we have just outlined, a natural algorithm is using the formulas~\eqref{caretodare} to convert~\eqref{care} into an equivalent~\eqref{dare} and solving it using~\eqref{sda-formulas}. This algorithm has been suggested by Chu, Fan and Lin~\cite{chu-care} as a doubling algorithm for CAREs. This algorithm inherits all the nice convergence properties of SDA for DAREs; in particular, among them, the fact that it also works (at reduced linear speed) on problems in which $A-GX_+$ has eigenvalues on the imaginary axis~\cite{Chi09-critical-doubling}.

While SDA works well in general, a delicate point is the choice of the shift value $\tau$. In principle almost every choice of $\tau$ works, since $\mathcal{H}-\tau I$ is singular only for at most $2n$ values of $\tau$, but in practice choosing the wrong value of $\tau$ may affect accuracy negatively. Dangers arise not from singularity of $\mathcal{H}-\tau I$ (which is actually harmless with a matrix pencil formulation), but from singularity in~\eqref{caretodare}, and also from taking $\tau$ too large or too small by orders of magnitude. A heuristic approach based on golden section search has been suggested~\cite{chu-care}.

In practice, one would prefer to avoid the Cayley transform or at least delay it as much as possible; this observation leads to another popular algorithm for CAREs. We start from the following observation.

\begin{lemma}
If $\mathcal{S} = c(\mathcal{H})$ (with a parameter $\tau \in \mathbb{R}$), then
\begin{equation} \label{squaringsignrelation}
	\mathcal{S}^2 = c\left(\frac{1}{2}\left(\mathcal{H} + \tau^2\mathcal{H}^{-1}\right)\right).	
\end{equation}
\end{lemma}
This identity can be verified directly, using the fact that rational functions of the same matrix $\mathcal{H}$ all commute with each other.

Applying this identity repeatedly, we get $\mathcal{S}^{2^k} = c(\mathcal{H}_k)$, where
\begin{align} \label{sign}
	\mathcal{H}_{k+1} &= \frac12 \left(\mathcal{H}_k + \tau^2\mathcal{H}_k^{-1}  \right), & \mathcal{H}_0 &= \mathcal{H}.
\end{align}
Hence one can hold off the Cayley transform and just compute the sequence $\mathcal{H}_k$ directly, starting from~\eqref{hamiltonian}. This constructs a sequence which represents implicitly $\mathcal{S}^{2^k}$. 

Constructing the matrices $\mathcal{H}_k$ is numerically much less troublesome than constructing explicitly $\mathcal{S}^{2^k}$ or its inverse $\mathcal{S}^{-2^k}$. Indeed, it is instructive to consider the behaviour of these iterations in a basis in which $\mathcal{H}$ is diagonal (when it exists). Let $\lambda$ be a generic diagonal entry (i.e., an eigenvalue) of $\mathcal{H}$. Then, $\mathcal{S}=c(\mathcal{H})$ has the corresponding eigenvalue $c(\lambda)$, and $\mathcal{S}^{2^k}$ has the eigenvalue $c(\lambda)^{2^k}$. If $\lambda \in \LHP$, then $\abs{c(\lambda)}<1$ (Lemma~\ref{lem:cayley}), and hence $c(\lambda)^{2^k} \to 0$ when $k\to\infty$. Similarly, if $\lambda$ is in the right half-plane, then $\abs{c(\lambda)}>1$ and $c(\lambda)^{2^k} \to \infty$. Thus $\mathcal{S}^{2^k}$ (as well as its inverse) has some eigenvalues that converge to zero, and some that diverge to infinity, as $k$ grows. This is one of the reasons why it is preferable to keep $\mathcal{S}$ in its factored form~\eqref{ssffact}. On the other hand, the eigenvalues of $\mathcal{H}_{k}$ converge to finite values $c^{-1}(0) = -\tau$ and $c^{-1}(\infty) = \tau$, so this computation suggests that the direct computation of $\mathcal{H}_k$ is feasible.

The \emph{sign function method}\cite{Rob80,denbea,GarL} to solve CAREs consists exactly in computing the iteration~\eqref{sign} up to convergence, obtaining a matrix $\mathcal{H}_\infty = \lim_{k\to\infty} \mathcal{H}_k$ that has numerically $n$ eigenvalues equal to $\tau\in \RHP$ and $n$ equal to $-\tau\in\LHP$, and then computing
\begin{align} \label{sign-final}
\operatorname{Im} \begin{bmatrix}
    U_1\\
    U_2
\end{bmatrix} &= \ker (\mathcal{H}_\infty + \tau I), & U_1,U_2& \in\mathbb{C}^{n\times n}, &  X_+ &= U_2U_1^{-1}.
\end{align}
The method takes its name from the fact that the limit matrix $\mathcal{H}_{\infty}$ (for $\tau=1$) is the so-called \emph{matrix sign function} of $\mathcal{H}$. We refer the reader to its analysis in Higham~\cite[Chapter~5]{higham-functions}, in which one clearly sees that one of the main ingredients is the formula~\ref{squaringsignrelation} relating the iteration to repeated squaring. 

Scaling is an important detail that deserves a discussion. Replacing $\mathcal{H}$ with a positive multiple of itself corresponds to multiplying each term of~\eqref{care} by a positive quantity; this operation does not change the solutions of the equation, nor the maximal / stabilizing properties of $X_+$. In SDA, scaling is limited to choosing the parameter of the initial Cayley transform, but in the sign method we have more freedom: we can take a different $\tau_k$ at each step of~\eqref{sign}. We remark that scaling for the sign method is usually presented in the literature in a slightly different form: one replaces~\eqref{sign} with 
\begin{equation} \label{sign2}
\mathcal{H}_{k+1} = \frac12\left((\tau_k^{-1} \mathcal{H}_k) + (\tau_k^{-1} \mathcal{H}_k)^{-1}\right).	
\end{equation}
The two forms are essentially equivalent, as they return iterates $\mathcal{H}_k$ that differ only by a multiplicative factor, which is then irrelevant in the final step~\eqref{sign-final}. Irrespective of formulation, the main result is that a judicious choice of scaling can speed up the convergence of~\eqref{sign} or~\eqref{sign2}. A cheap and effective choice of scaling, \emph{determinantal scaling}, $\tau_k = (\det \mathcal{H}_k)^{\frac1{n}}$ has been suggested by Byers~\cite{Bye87}. Other related choices of scaling and their performances have been discussed by Higham~\cite[Chapter~5]{higham-functions} and Kenney and Laub~\cite{KenL}. The general message is that scaling has a great impact in the first steps of the iteration, when it can greatly improve convergence, but once the residual starts to decrease its effect in the later steps becomes negligible.

Scaling also has an impact on stability; the stability of the sign iteration as a method to compute invariant subspaces (and hence ultimately Riccati solutions) has been studied by Bai and Demmel~\cite{BaiD98} and Byers, He and Mehrmann~\cite{ByeHM97}. The two interesting messages are that (expectedly) the sign function method suffers when $\mathcal{H}$ is ill-conditioned, but that (unexpectedly) the invariant subspaces extracted from $\mathcal{H}_\infty$ has better stability properties than $\mathcal{H}_\infty$ itself. A version of the sign iteration that uses matrix pencils to reduce the impact of these inversions have been suggested by Benner and Byers~\cite{BenB06}.

Another useful computational detail is that one can rewrite the sign function method~\eqref{sign} as
\begin{align*}
\mathcal{M}_{k+1} &= \frac12 (\mathcal{M}_k + \tau^2 J \mathcal{M}_k^{-1}J), &  \mathcal{M}_k = \mathcal{H}_k J,
\end{align*}
which is cheaper because one can take advantage of the fact that the matrices $\mathcal{M}_k$ are Hermitian~\cite{Bye87}. Indeed, it is a general observation that most of the matrix algebra operations needed in doubling-type algorithms can be reduced to operations on symmetric/Hermitian matrices; see for instance also~\eqref{caretodare}.

\subsection{Remarks}

The formulation in the sign iteration allows one to introduce some form of per-iteration scaling in the setting of a doubling-type algorithm. It would be interesting to see if this scaling can be transferred to the SDA setting, and which computational advantage it brings. Note that, in view of~\eqref{squaringsignrelation}, scaling the sign iteration is equivalent to changing the parameter $\tau$ in the Cayley transform. So SDA does incorporate a form of scaling, but only at the first iteration, when one chooses $\tau$.

In general, it is unclear if scaling after the first iteration produces major gains in convergence speed. It would be appealing to try and study this kind of scaling with the tools of polynomial and rational approximation, like it has been done in more details for non-doubling algorithms, with the aim of deriving optimal choices for the parameters $\tau$ and $\sigma_k$.

There is another classical iterative algorithm to solve algebraic Riccati equations (both in discrete and continuous time), and it is Newton's method. For the simpler case of CAREs, Newton's method~\cite{Kle68} consists in determining $X_{k+1}$ by solving at each step the Lyapunov equation
\begin{equation} \label{newton-step}
	(A-GX_k)^*(X_{k+1}-X_k) + (X_{k+1}-X_k)(A-GX_k) = -(Q+A^*X_k+X_kA-X_k G X_k)	
\end{equation}
or the equivalent one
\[
(A-GX_k)^*X_{k+1} + X_{k+1}(A-GX_k) = -Q - X_k G X_k.
\]
A line search procedure, which improves convergence speed in practice, has been introduced by Benner and Byers~\cite{BenB98-linesearch}. The method can be used, in particular, for large and sparse equations in conjunction with low-rank ADI~\cite{BenLP08}.

The reader may wonder if there is an explicit relation between doubling algorithms and Newton-type algorithms, considering especially that both exhibit quadratic convergence (which, moreover, in both cases degrades to linear with rate $1/2$ if $A-GX_+$ has purely imaginary eigenvalues~\cite{GuoL98}). The answer, unfortunately, seems to be no. An argument that suggests that the two iterations are genuinely different is that the iterates produced by Newton's method approach $X_+$ from \emph{above}~\cite{Kle68} (i.e., $X_1 \succeq X_2 \succeq \dots \succeq X_{k}\succeq X_{k+1} \succeq \dots \succeq X_+$), not from \emph{below} like the iterates $Q_k$ of SDA in~\eqref{sda-formulasQ}.

Some more recent algorithms for large and sparse CAREs essentially merge the Newton step~\eqref{newton-step} and the ADI iteration~\eqref{adi-step} into a single iteration~\cite{LinSim,Sim16,BenBKS}. It is again unclear whether there is an explicit relation between these two families of methods.

An interesting question is what is the `non-doubling' analogue of the sign method and of SDA. One can convert the CARE to discrete-time using~\eqref{caretodare} and formulate~\eqref{dare-iterative}, but to our knowledge this method does not have a more appealing presentation in terms of a simple iterative method for~\eqref{care}, like it has in all the other discrete-time examples.

Another `philosophical' observation is that the sign function method does not avoid a Cayley-type transformation; it merely pushes it back to the very last step~\eqref{sign-final}, where the sub-expression $\mathcal{H}+ \tau I$ appears; this operation takes the role of a discretizing transformation that maps the eigenvalue $-\tau$ into a value inside a given circle and the eigenvalue $\tau$ into one outside. A discretizing transformation of some sort seems inevitable in this family of algorithms, although delaying it until the very last step seems beneficial for accuracy, because at that point we have complete control of the location of eigenvalues.

\section{Unilateral equations and NMEs} \label{sec:nme}

We end our discussion of the family of Riccati-type equations with a pair of oft-neglected cousins, and present them with an application that shows clearly the relationship between them. Consider the matrix Laurent polynomial
\begin{align} \label{laurpol}
P(z) &= Az^{-1} + Q + A^* z, & Q &= Q^* \succ 0, & A,Q &\in \mathbb{C}^{n\times n}.
\end{align}
The problem of \emph{spectral factorization} (of quadratic matrix polynomials) consists in determining a factorization
\begin{align} \label{specfact}
P(z) &= (z Y^* - I)X(z^{-1}Y - I), & X&=X^* \succ 0, & X,Y&\in\mathbb{C}^{n\times n},
\end{align}
such that $\rho(Y) \leq 1$. In particular, the left factor is invertible for $\abs{z}<1$, and the right factor is invertible for $\abs{z}>1$.

Equating coefficients in~\eqref{laurpol} and~\eqref{specfact} gives $-XY=A$, $Q = X + Y^*XY$. We can eliminate one among $X$ and $Y$ from this system of two equations, getting two equations with a single unknown each
\begin{align}
0 &= A + QY + A^*Y^2, \label{uqme} \\
Q &= X + A^*X^{-1}A. \label{nme}
\end{align}
The first one~\eqref{uqme} is called \emph{unilateral quadratic matrix equation}~\cite{BinILM06}, while the second one~\eqref{nme} is known with the (rather undescriptive) name of \emph{nonlinear matrix equation} (NME)~\cite{GuoL-nano1,GuoL-nano2,doublingbook}.

While~\eqref{uqme} looks more appealing at first, as it reveals direct ties with the palindromic quadratic eigenvalue problem~\cite{GuoL-nano1,GuoL-nano2,M4}, it is in fact~\eqref{nme} that reveals more structure: for instance, \eqref{nme}~has Hermitian solutions (see below), while the structure in the solutions of~\eqref{uqme} is much less apparent.

\subsection{Solution properties}

It follows from~\eqref{specfact} that $P(\lambda) \succeq 0$ for each $\lambda$ that belongs to the unit circle (hence $\lambda^{-1} = \bar{\lambda}$), so this is a necessary condition for the solvability of this problem. It can be proved that it is sufficient, too, and that a maximal / stabilizing solution exists.

\begin{theorem}\cite[Theorem~2.2]{EngRR} \label{thm:nmesol} Assume that $P(z)$ is regular and $P(\lambda) \succeq 0$ for each $\lambda$ on the unit circle. Then, \eqref{nme} has a (unique) solution $X_+$ such that
\begin{enumerate}
	\item[(1)] $X_+ = X_+^* \succ 0$;
	\item[(2)] $X_+ \succeq X$ for any other Hermitian solution $X$;
	\item[(3)] $\rho(Y) = \rho(-X_+^{-1}A) \leq 1$
\end{enumerate}
If, in addition, $P(\lambda) \succ 0$ for each $\lambda$ on the unit circle, then $\rho(-X_+^{-1}A) < 1$.
\end{theorem}

Once again, we can rewrite~\eqref{nme} as an invariant subspace problem.
\begin{align}
\begin{bmatrix}
    A & 0\\
    -Q & I
\end{bmatrix}
\begin{bmatrix}
    I\\X
\end{bmatrix}
&=
\begin{bmatrix}
    0 & -I\\
    A^* & 0
\end{bmatrix}
\begin{bmatrix}
    I\\X
\end{bmatrix}
Y, & Y &= -X^{-1}A.
\end{align}

We assume again that $A$ is invertible to avoid technicalities with matrix pencils. The matrix
\begin{equation}
	\mathcal{S} = \begin{bmatrix}
	    0 & -I\\
	    A^* & 0
	\end{bmatrix}^{-1}
	\begin{bmatrix}
	    A & 0\\
	    -Q & I
	\end{bmatrix}
\end{equation}
is symplectic, and so is the slightly more general form
\begin{equation} \label{symplectic2}
	\begin{bmatrix}
	    G & -I\\
	    A^* & 0
	\end{bmatrix}^{-1}
	\begin{bmatrix}
	    A & 0\\
	    -Q & I
	\end{bmatrix}.
\end{equation}
\begin{lemma}
\begin{enumerate}
	\item[(1)] A matrix in the form~\eqref{symplectic2} is symplectic if and only if $G=G^*, Q=Q^*$, and the two blocks called $A,A^*$ in~\eqref{symplectic} are one the conjugate transpose of the other.
	\item[(2)] If $\lambda$ is an eigenvalue of a symplectic matrix with right eigenvector $v$, then $\overline{\lambda}^{-1}$ is an eigenvalue of the same matrix with left eigenvector $v^*J$.
	\item[(3)] If the hypotheses of Theorem~\ref{thm:nmesol} hold (including the strict positivity one in the end), then the $2n$ eigenvalues of $\mathcal{S}$ are (counting multiplicities) the $n$ eigenvalues $\lambda_1,\lambda_2,\dots,\lambda_n$ of $-X_+^{-1}A$ inside the unit circle, and the $n$ eigenvalues $\overline{\lambda_i}^{-1}$, $i=1,2,\dots,n$ outside the unit circle.
\end{enumerate}
\end{lemma}
The symplectic structure behind this equation is the same one as the DARE, and indeed Part~2 of this lemma is identical to Part~2 of Lemma~\ref{symplemma}. Indeed, Engwerda, Ran and Rijkeboer~\cite[Section~7]{EngRR} note that~\eqref{nme} can be reduced to a DARE, although it is one that does not fall inside our framework since it has $G \preceq 0$.

\subsection{Algorithms}

The formulation~\eqref{nme} suggests immediately the iterative algorithm
\begin{equation} \label{nme-iterative}
	X_{k+1} = Q - A^* X_k^{-1}A.
\end{equation}
Clearly we cannot start this iteration from $0$, so we take $X_1 = Q$ instead. An interesting interpretation of this algorithm is as iterated Schur complements of block Toeplitz tridiagonal matrices. The Schur complement of the $(1,1)$ block of the tridiagonal matrix
\[
\underbrace{
\begin{bmatrix}
    X_{k} & A^*\\
    A & Q & A^*\\
    & A & Q & \ddots\\
    & & \ddots & \ddots & A^*\\
    & & & A & Q
\end{bmatrix}}_{\text{$h$ blocks}},
\]
is 
\[
\underbrace{
\begin{bmatrix}
    X_{k+1} & A^*\\
    A & Q & A^*\\
    & A & Q & \ddots\\
    & & \ddots & \ddots & A^*\\
    & & & A & Q
\end{bmatrix}}_{\text{$h-1$ blocks}}.
\]
Hence the whole iteration can be interpreted as constructing successive Schur complements of the tridiagonal matrix
\begin{equation} \label{tridiag}
\mathcal{Q}_{m} :=
\underbrace{\begin{bmatrix}
    Q & A^*\\
    A & Q & A^*\\
    & A & Q & \ddots\\
    & & \ddots & \ddots & A^*\\
    & & & A & Q
\end{bmatrix}}_{\text{$m$ blocks}}.
\end{equation}
It can be seen that $\mathcal{Q}_m$ is positive semidefinite, under the assumptions of Theorem~\ref{thm:nmesol}: a quick sketch of a proof is as follows. The matrix $\mathcal{Q}_m$ is a submatrix of
\[
\begin{bmatrix}
    Q & A^* & & & A\\
    A & Q & A^*\\
    & A & Q & \ddots\\
    & & \ddots & \ddots & A^*\\
    A^*& & & A & Q
\end{bmatrix}
= (\Phi \otimes I) \begin{bmatrix}
    P(1)\\
    & P(\zeta)\\
    & & P(\zeta^2)\\
    & & & \ddots\\
    & & & & P(\zeta^{-1})
\end{bmatrix} (\Phi \otimes I)^{-1},
\]
which the equation shows to be similar (using the Fourier matrix $\Phi$ and properties of Fourier transforms) to a block diagonal matrix that contains $P(z)$ from~\eqref{laurpol} evaluated in the roots of unity $1,\zeta,\zeta^2,\dots,\zeta^{-1}$. 

Hence, in particular, all the $X_{k}$ are positive semidefinite. One can further show that $Q = X_0 \succeq X_1 \succeq X_2 \succeq \dots  \succeq X_k \succeq \dots$. The sequence $X_k$ is monotonic and bounded from below, hence it converges, and one can show that its limit is $X_+$~\cite[Section~4]{EngRR} (to do this, verify the property in Point~(2) of Theorem~\ref{thm:nmesol} by proving that $X_k \succeq X$ at each step of the iteration).

A doubling variant of~\eqref{nme-iterative} can be constructed starting from this Schur complement interpretation. The Schur complement of the submatrix formed by the odd-numbered blocks $(1,3,5,\dots,2m-1)$ of
\[
\underbrace{
\begin{bmatrix}
    U_{k} & A_k^*\\
    A_k & U_k & A_k^*\\
    & A_k & \ddots & \ddots\\
    & & \ddots & U_k & A_k^*\\
    & & & A_k & Q_k
\end{bmatrix}}_{\text{$2m$ blocks}},
\]
is
\[
\underbrace{
\begin{bmatrix}
    U_{k+1} & A_{k+1}^*\\
    A_{k+1} & U_{k+1} & A_{k+1}^*\\
    & A_{k+1} & \ddots & \ddots\\
    & & \ddots & U_{k+1} & A_{k+1}^*\\
    & & & A_{k+1} & Q_{k+1}
\end{bmatrix}}_{\text{$m$ blocks}},
\]
with
\begin{subequations} \label{cr}
\begin{align}
A_{k+1} &= -A_kU_k^{-1}A_k,\\
Q_{k+1} &= Q_k - A_k^* U_k^{-1} A_k,\\
U_{k+1} &= U_k - A_k^* U_k^{-1}A_k - A_k U_k^{-1}A_k^*.
\end{align}
\end{subequations}
We can construct the Schur complement of the first $2^k-1$ blocks of $\mathcal{Q}_{2^k}$ in two different ways: either we make $2^k-1$ iterations of~\eqref{nme-iterative}, resulting in $X_{2^k}$, or we make $k$ iterations of~\eqref{cr}, starting from $A_0=A, Q_0=U_0=Q$, resulting in $Q_k$. This shows that $Q_k = X_{2^k}$.

This peculiar way to take Schur complements of Toeplitz tridiagonal matrices was introduced by Buzbee, Golub and Nielson~\cite{BuzGN} to solve certain differential equations, and then later applied to matrix equations similar to~\eqref{uqme} and~\eqref{nme} by Bini, Gemignani, and Meini~\cite{BM,BGM,Mei}. The iteration~\eqref{cr} is known as \emph{cyclic reduction}.

One can derive the same iteration from repeated squaring, in the same way as we obtained SDA as a modified subspace iteration~\cite{LinX06}. We seek formulas to update a factorization of the kind
\[
\mathcal{S}^{-2^k} = \begin{bmatrix}
    A_k & 0\\
    -Q_k & I
\end{bmatrix}^{-1}
\begin{bmatrix}
    G_k & -I\\
    A_k^* & 0
\end{bmatrix}.
\]
To do this, we write (analogously to~\eqref{sdaderiv})
\[
\mathcal{S}^{-2^{k+1}} = \mathcal{S}^{-2^k} \mathcal{S}^{-2^k} =
\begin{bmatrix}
    A_k & 0\\
    -Q_k & I
\end{bmatrix}^{-1}
\left(
\begin{bmatrix}
    G_k & -I\\
    A_k^* & 0
\end{bmatrix}
\begin{bmatrix}
    A_k & 0\\
    -Q_k & I
\end{bmatrix}^{-1}
\right)
\begin{bmatrix}
    G_k & -I\\
    A_k^* & 0
\end{bmatrix}
\]
and use Lemma~\ref{lem:bmf} (with $[M_1\, M_2] = I_{2n}$) to find a factorization in the form~\eqref{bmf-factorization} of the term in parentheses, which then combines with the outer terms to produce the sought decomposition. The resulting formulas are
\begin{subequations} \label{sda2}
\begin{align}
A_{k+1} &= -A_k(Q_k-G_k)^{-1}A_k,\\
Q_{k+1} &= Q_k - A_k^* (Q_k-G_k)^{-1} A_k,\\
G_{k+1} &= G_k + A_k (Q_k-G_k)^{-1} A_k^*,
\end{align}
\end{subequations}
and one sees that they coincide with~\eqref{cr}, after setting $U_k = Q_k - G_k$. With an argument analogous to the one in Section~\ref{sec:dare}, one sees that
\[
\mathcal{S}^{-2^k}\begin{bmatrix}
    0 \\ -I
\end{bmatrix} = \begin{bmatrix}
    I\\
    Q_k
\end{bmatrix},
\]
thus $\begin{bsmallmatrix}
    I\\Q_k
\end{bsmallmatrix}$ converges to a basis of the invariant subspace associated to the eigenvalues of $\mathcal{S}$ inside the unit circle.

This formulation~\eqref{sda2} is known as \emph{SDA-II}~\cite{LinX06,Chu-trains}. 

\subsection{Remarks}

Even though we have mentioned spectral factorization only here, it can be formulated for more complicated matrix functions also in the context of DAREs and CAREs; in fact, it is a classical topic, and another facet of the multiple connections between matrix equations and control theory~\cite{bart1,bart2,SayK}.

The interpretation as Schur complement is a powerful trick, which reveals a greater picture in this family of methods. It may possibly be used to understand more about the stability of these methods, since Schur complementation and Gaussian elimination on symmetric positive definite matrices is a well understood topic from the numerical point of view.

Many authors have studied variants of~\eqref{nme}. Typically, one replaces the nonlinear term with various functions of the form $A^* f(X) A$, or adds more nonlinear terms. In the modified versions, it is often possible to prove convergence of the fixed-point algorithm with arguments of monotonicity, and prove the existence of a solution under some assumptions. However, after any nontrivial modification the connection with invariant subspaces is lost. This fact, coupled with lack of applications, makes these variants much less interesting than the original equation, in the eyes of the author.

% \subsection{Spectral division}

% [TODO: Malyshev and Bai-Demmel-Gu?]

% [TODO: Nakatsukasa - polar?]

\section{Nonsymmetric variants in applied probability}

Many of the equations treated here have nonsymmetric variants which appear naturally in queuing theory, a sub-field of applied probability. In the analysis of \emph{quasi-birth-death models}~\cite{LR-book,BLM-book}, one encounters equations of the form
\begin{align} \label{qbd}
0 &= A + QY + BY^2, & A, B, Q, Y \in \mathbb{R}^{n\times n}, 
\end{align}
where $A,B \geq 0$ (we use the notation $M\geq N$ to denote that a matrix $M$ is entrywise larger than $N$, i.e., $M_{ij} \geq N_{ij}$ for all $i,j$), and the matrix $-Q$ is an M-matrix, i.e., $Q_{ij} \geq 0$ for $i\neq j$ and $\Lambda(Q) \subset \overline{\LHP}$. These equations have a solution $Y \geq 0$ which has a natural probabilistic interpretation. The solution $X$ to $X = Q - BX^{-1}A$ and the solution of the associated dual equation $0 = Z^2A + ZQ+B$ also appear naturally and have a related probabilistic meaning~\cite[Chapter~6]{LR-book}\cite[Section~5.6]{BLM-book}.

Similarly, the equation
\begin{align}  \label{nare}
Q + BX+XA - XGX &= 0, & Q,X & \in \mathbb{R}^{m\times n}, A \in \mathbb{R}^{n\times n}, B \in \mathbb{R}^{m\times m}, G\in\mathbb{R}^{n\times m}.
\end{align}
appears in the study of so-called \emph{fluid queues}, or \emph{stochastic flow models}~\cite{roger94,kk95,soares05}. The matrices $A,B$ are M-matrices, while $G,-Q \geq 0$. 
One can formulate nonsymmetric analogues of basic matrix iterations and doubling algorithms. Unfortunately, the theory does not translate perfectly to this setting, due to the sign differences between the two cases: in the symmetric equations $G,Q\succeq 0$, while in the nonsymmetric case $G,-Q \geq 0$. Due to this asymmetry, the signs in the two cases do not match, and one needs to formulate different arguments. For instance, in the symmetric case one proves that the inverses that appear in~\eqref{sda-formulas} exist because $G_k \succeq 0$, $Q_k \succeq 0$; while in its nonsymmetric analogue $G_k,-Q_k \geq 0$, and one proves that $I+G_kQ_k$ and $I+Q_kG_k$ are M-matrices to show that those inverses exist.

The equation~\eqref{dare} does not appear to have an immediate analogue in queuing theory, but this fact seems just an accident, since some of the results that involve~\eqref{nare} could have been formulated with an equivalent equation resembling more~\eqref{dare} than~\eqref{care} instead. There is a distinction between discrete-time and continuous-time models also in applied probability, but in many cases it does not affect directly the shape of the equations; for instance~\eqref{qbd} takes the same form for discrete- and continuous-time QBDs. The role of discretizing transformations such as Cayley transforms in this context has been studied by Bini, Meini, and Poloni~\cite{BMP-transforming}.

For reasons of space, we cannot give here a complete treatment of these nonsymmetric variants. Huang, Li and Lin~\cite{doublingbook} in their book enter into more detail about the doubling algorithms for these equations, but a great part of the theory (including existence results and probabilistic interpretations for the iterates of various numerical methods) is unfortunately available only in the queuing theory literature, strictly entangled with its applications.

An interesting remark is that the M-matrix structure allows one to construct stability proofs more easily. Conditioning and stability results for these equations have been studied by some authors~\cite{XXL1,XXL2,NguP15,XL17,CheLM19}, relying heavily on the sign and M-matrix structure. The forward stability proof in Nguyen and~Poloni~\cite{NguP15} is, to date, one of the very few complete stability proofs for a doubling-type algorithm.

\section{Conclusions}

In this paper, we presented from a consistent point of view doubling algorithms for symmetric Riccati-type equations, relating them to the basic iterations of which they are a `squaring' variant. We have included various algorithms that belong to the same family but have appeared independently, such as the sign iteration and cyclic reduction. We have outlined relations between doubling algorithms, the subspace iteration, ADI-type and Krylov subspace methods, and Schur complementation of tridiagonal block Toeplitz matrices. This theory, in turn, forms only a small portion of the far larger topic of numerical algorithms for Riccati-type equations and control theory. This field of research is an incredibly vast one, spanning at least six decades of literature and various communities between engineering and mathematics, so we have surely omitted or forgotten many relevant contributions; we apologize with the missing authors.

We hope that the reader can benefit from our paper by both gaining theoretical insight, and having available some numerical algorithms for these equations. Indeed, with respect to many competitors, doubling-based algorithms have the advantage that they reduce to the simple coupled matrix iterations~\eqref{sda-formulas} or~\eqref{cr}, which are easy to code and fast to run in many computational environments.
 
Another interesting remark that was suggested by a referee is that some recent lines of research consider this family of matrix equations under different types of data sparsity than low-rank:  for instance, Palitta and Simoncini~\cite{PaliS} consider banded data, and Kressner, Massei and Robol~\cite{KreMR} and Massei, Palitta and Robol~\cite{MasPR} consider semi-separable (low-rank off-diagonal blocks) and hierarchically semiseparable structures. Much earlier, Grasedyck, Hackbusch and Khoromskij~\cite{GraHK} considered using hierarchical matrices to solve Riccati equations. All these structures are (at least up to a degree) preserved by the operations involved in doubling methods~\cite{HierBook,XiaCGL}. These novel techniques may open up new lines of research for doubling-type algorithms.

% [TODO: Benner-Byers pencil arithmetic] 

\bibliographystyle{abbrvurl}
\bibliography{paper}
\end{document}